\documentclass[preprint,12pt]{elsarticle}

\usepackage{amssymb}
\usepackage{amsmath,amssymb,amsfonts}
\usepackage{amsthm}
\usepackage{array}
\usepackage{algorithm}
\usepackage{algorithmic}
\usepackage{booktabs}
\usepackage{color}
\usepackage{dsfont}
\usepackage{graphicx}
\usepackage{hyperref}
\usepackage{makecell}
\usepackage{setspace}

\usepackage{subfigure}

\usepackage{stfloats}
\usepackage{txfonts}
\usepackage{times}
\usepackage{textcomp}
\usepackage{url}
\usepackage{verbatim}

\journal{Applied Mathematics and Computation}

\begin{document}

\begin{frontmatter}



\title{Two kinds of parametric piecewise rational  interpolation kernels  for image magnification}


\author[label1]{Bing Guo}
\affiliation[label1]{organization={School of Mathematics and Statistics, Jishou University}, 
	addressline={No.120, Renmin South Road},
	city={Jishou},
	postcode={416000},
	state={Hunan Province},
	country={China}}

\author[label2]{Wanfeng Qi}

\affiliation[label2]{organization={School of Mathematics, Liaoning Normal University}, 
            addressline={No.850, Huanghe Road},
            city={Dalian},
            postcode={116029},
            state={Liaoning Province},
            country={China}}

\begin{abstract}
We study the constructions of piecewise rational interpolation kernels that are supported on the interval $[-2,2]$,
and present one novel rational cubic/linear and  five quartic/linear interpolation kernels.
All proposed kernels are symmetric, $C^1$ continuous, and possess certain degrees of approximation order.
The proposed quartic/linear interpolation kernels include the cubic and the cubic/linear interpolation kernel as special cases.
Our numerical results show that one of the quartic/linear interpolation kernels can outperform the cubic interpolation kernel in terms of PSNR, SSIM, and FSIM.

\end{abstract}




\end{frontmatter}


\section{Introduction}

Interpolation  provides a means of estimating   function values at intermediate points from discrete sampling points,
and  is widely used in image processing.
It plays a crucial role in resampling and reconstructing images, but image quality often deteriorates during image processing, primarily depending on the properties of the interpolation kernel.

Interpolation functions with various interpolation kernels
have been introduced, typically expressed in the form
\begin{eqnarray*}
	g(t)=\sum_{i \in \mathbb{Z}} c_{i}\phi(t-t_{i}),
\end{eqnarray*}
where $t_{i}$ represents the time sampling points, $\phi(t)$
denotes the interpolation kernel, and ${c_{i}}$ represents given discrete samples.
Sinc function  is the well-known  ideal interpolation kernel which is defined by
\begin{equation*}
	\text{sinc}(t)=\left\{
	\begin{aligned}
		&\frac{\sin (\pi t)}{\pi t}, ~~~&& t\neq 0,\\
		&1, ~~~&&t=0.
	\end{aligned}
	\right.
\end{equation*}
However, the ideal interpolation kernel is seldom used in practice due to its infinite impulse response and slow decay at infinity.
To address this issue, various kernels of finite support have been proposed.


Several interpolation kernels with parameters are widely used in image magnification. The rectangular kernel  \cite{RifimanM74, Berstein76}  and triangular kernel \cite{LehmannGS99} are the simplest kernels of finite support.
The cubic spline convolution kernel, originally developed for Landsat digital image reconstruction, is another well-known kernel of finite support.
It offers high reconstruction accuracy with relatively low computational complexity.
A family of  functions constructed in a manner similar to spline has received widespread attention and is considered as an ideal replacement for the sinc kernel.
Comprehensive details about these kernels are available in the works \cite{LehmannGS99, ParkerKT83, maeland88, meijeringZV99}.
The following section will provide specific details on the three classical kernels of finite support mentioned above, which are widely utilized in computer-assisted surgery, medical photographs, clinical images, remote sensed data, and other related fields \cite{Simon75, hanB00, shiR06, tsaiLT11, banerjeeRS12,ZhouSD12, hongWT13}.
Polynomial and spline interpolation are frequently used in nonuniform sampling cases as an alternative to convolution methods in uniform sampling cases.
Zou et al. \cite{zou2020new} developed a Newton-Type polynomial interpolation with adjustable parameters to modify the value of the interpolation polynomial. Parametric rational spline interpolation has been explored for shape preservation, including positivity and monotonicity.
Delbourgo and Gregory \cite{delbourgo1985shape} investigated a parametric $C^1$ rational cubic function to address the issue of shape preserving interpolation. Han \cite{han2008convexity}  proposed a parametric piecewise rational quartic interpolant with a suitable selection of parameters to maintain the local convexity of the given data.
Additionally, Han \cite{han2015shape} analyzed another type of parametric shape preserving piecewise rational interpolant, which has a quartic numerator and quadratic denominator, and provided an explicit representation of the parameters.
Sarfraz et al. \cite{sarfraz1997preserving}  introduced a $C^1$  monotonically shape preserving scheme by the data dependent choice of shape parameters.


Piecewise rational interpolation kernels show promise as alternatives to classical polynomial interpolation kernels in uniform sampling cases.
However, constructing meaningful rational interpolation kernels poses a challenge in ensuring they satisfy the partition property of unity.
To the best of our knowledge, the only existing rational kernel function is Hu and Tan's quadratic/quadratic rational interpolation kernels $S_{2/2}$ \cite{HuT06}. It is based on $4$th Thiele-type interpolating continued fractions of sinc function and can only be constructed in segments. In each segment, it interpolates four points of the sinc function and expresses them in the form of a continued fraction, combined with one derivative interpolation condition. Although it has one parameter, it does not possess the partition of unity property.
And for each parameter, it should be recalculated to get the final  explicit expression.
In contrast, we use a direct approach to convert the partition of unity requirement into solving a system of nonlinear equations, resulting in new rational interpolation kernel functions with explicit expressions and partition of unity property.
Specifically, we first express the partition of unity property as an expression in interval $[0,1]$, then convert it into a rational function that equals 1.
Next, we require that the numerator equals to the denominator of the rational function, resulting in a complicated system of nonlinear equations.
We solve these equations and substitute the solutions back into the original form.
Finally, we rigorously check the final form for any potential singularities in the denominators.

The following are our main contributions:
\begin{itemize}
	\item  We propose one new  rational cubic/linear interpolation kernel
	and five  quartic/linear interpolation kernels.
	These kernels possess explicit forms and simultaneously satisfy the partition of unity property and $C^1$ continuity,
	while Hu and Tan's $S_{2/2}$ \cite{HuT06}
	only provide an explicit form for a specific parameter, 
	which does not satisfy the partition of unity property. 
	When the parameter takes other values, it needs to be recalculated.
	\item The proposed quartic/linear interpolation kernels, which include the well-known cubic interpolation kernel as a special case, have free parameters that provide higher flexibility in various applications.
    Numerical experiments demonstrate that one quartic/linear interpolation kernels with fixed parameters outperform the cubic interpolation kernel in image magnification.
\end{itemize}


The paper is structured as follows: Section \ref{section 2} provides an overview of three classical interpolation kernels, namely the nearest neighbor, linear and cubic interpolation kernels. Sections \ref{sec:Rational cubic convolution interpolation} and \ref{section:Rational quartic/linear interpolation kernel $S_{4/1}(t)$} introduce the rational cubic/linear kernel and quartic/linear interpolation kernels, respectively. Section \ref{sec:experiments} presents the results of numerical experiments conducted on two quartic/linear interpolation kernels.





\section{Classical	interpolation kernels}
\label{section 2}

The nearest neighbor interpolation kernel is given by
\begin{equation}
	\label{eqn:nearest neighbor kernel}
	S_0(t)=\left\{
	\begin{aligned}
		&1, ~~~&& 0\leq |t|<0.5,\\
		&0, ~~~&& |t|\geq 0.5.
	\end{aligned}
	\right.
\end{equation}

The linear interpolation kernel $S_1(t)$ is  represented by the triangular function
\begin{equation}
	\label{eqn: linear kernel}
	S_1(t)=\left\{
	\begin{aligned}
		&1-|t|, ~~~&& 0\leq |t|<1,\\
		&0, ~~~&& |t|\geq 1,
	\end{aligned}
	\right.
\end{equation}
which can be viewed as a linear approximation of the sinc function.

The initial cubic interpolation kernel $S_3(t)$ with eight undetermined parameters is
\begin{equation*}
	S_3(t)=\left\{
	\begin{aligned}
		&a_{00}+a_{01}|t|+a_{02}|t|^{2}+a_{03}|t|^{3}, ~~~&& 0\leq |t|<1,\\
		&a_{10}+a_{11}|t|+a_{12}|t|^{2}+a_{13}|t|^{3}, ~~~&& 1\leq |t|<2,\\
		&0, ~~~&&  |t|\geq 2.
	\end{aligned}
	\right.
	\label{eqn:cubic with 9 parameter}
\end{equation*}

There are several indicators to assess the quality of interpolation kernels, one of which is the degree of  approximation order.
A kernel  $\phi(t)$ satisfying $\int_{-\infty}^{\infty} \phi(t)dt\neq 0$ is said to have approximation order $L$
if  for  any polynomial $P(t)$ of degree less than or equal to $L-1$, there exists $C_P\in \mathbb{R}$ such that
\begin{equation}
	\label{eqn:approximation order definition}
	\sum_i P(t-i)\phi(t-i)=Cp
\end{equation}
holds  almost everywhere \cite{blu2003complete}.
This condition is closely related to the well-known Strang-Fix condition \cite{strang2011fourier}, which is an important property for evaluating various interpolation kernels.
If the interpolation kernel function has a rational form, equation (\ref{eqn:approximation order definition}) usually leads to a complicated system of nonlinear equations that are often difficult to solve. This poses a challenge to constructing interpolation kernel functions in a rational form that possess high-order approximation order.
If $L=1$, it is well-known as  the partition of unity property
\begin{equation}
	\label{eqn:partition of unity}
	\sum_{i=-\infty}^\infty \phi(t-i)=1,~~~~~t\in \mathbb{R}.
\end{equation}
Both the nearest neighbor interpolation kernel $S_0(t)$ and the linear interpolation kernel $S_1(t)$ satisfy this property.

If the cubic interpolation kernel $S_3(t)$ satisfies the partition of unity (\ref{eqn:partition of unity}) and $C^{1}$ continuous condition,
then it has only one free parameter and becomes
	\begin{equation}
		\label{eqn:cubic interpolation kernel}
		S_3(t)=\left\{
		\begin{aligned}
			&(1-|t|)(1+|t|+(1+a_{02})|t|^2), & 0\leq |t|<1,\\
			&-(3+a_{02})(2-|t|)^2(1-|t|),\!& 1\leq |t|<2,\\
			&0, &  |t|\geq 2.
		\end{aligned}
		\right.
	\end{equation}
Several choices for the parameter $a_{02}$ have been analyzed. Keys \cite{Keys81} recommended the value of $a_{02}=-5/2$ for the cubic interpolation kernel $S_3(t)$ to achieve an approximation order of 3.
Park and Schowengerdt \cite{ParkS83} have extensively discussed various choices of $a_{02}$, including $a_{02}=-3,-5/2,-\frac{7}{3},-2$.
When $a_{02}=-2$, we have
$$
S_3^{'}(1)=\text{sinc}'(1)=-1.
$$

\section{Rational cubic/linear interpolation kernel \(S_{3/1}(t)\)}
\label{sec:Rational cubic convolution interpolation}
The initial symmetric rational cubic/linear interpolation kernel $S_{3/1}(t)$, which has ten unknown parameters, is defined as
\begin{equation}
	\label{eqn:cubic/linear initial form}
	S_{3/1}(t)=
	\left\{\begin{array}{ll}
		\frac{a_{00}+a_{01}|t|+a_{02}|t|^{2}  +a_{03}|t|^{3}}{1 +b_{01}|t|}, & 0 \leq |t|<1,\\
		\frac{a_{10}+a_{11}|t| +a_{12}|t|^{2}  +a_{13}|t|^{3}}{1+b_{11}|t|}, & 1 \leq |t|<2,\\
		0,& |t|\geq 2.
	\end{array}\right.
\end{equation}
It has a cubic numerator and a linear denominator.
Similar to $S_3(t)$, $S_{3/1}(t)$ should interpolate the sinc function at $t=0, 1, 2$, i.e.,
$S_{3/1}(0)=\text{sinc}(0)=1,
S_{3/1}(1)=\lim_{t\to1^{-}}S_{3/1}(t)=\text{sinc}(1)=0,
S_{3/1}(2)=\lim_{t\to2^{-}}S_{3/1}(t)=\text{sinc}(2)=0$,
which are equivalent to
\begin{equation*}
	\left\{
	\begin{aligned}
		&a_{00}= 1,\\
		& \frac{a_{00}+a_{01}+a_{02}+a_{03}}{1 + b_{01}} = 0,\\
		&\frac{a_{10}+a_{11}+a_{12}+a_{13}}{
			1 + b_{11}}= 0,\\
		&\frac{	a_{10}+2a_{11}+4a_{12}+8a_{13}}{1 +2b_{11}} =0.
	\end{aligned}
	\right.
\end{equation*}
Under the conditions $b_{01}\neq -1, b_{11}\neq -1, b_{11}\neq -1/2$, we can get
$a_{00}=1, a_{02}=-1-a_{01}-a_{03}, a_{11}=2a_{13}-3a_{10}/2, a_{12}=-3a_{13}+a_{10}/2$, and
the following form
\begin{equation}
	\label{eqn:cubic/linear interpolation form}
	S_{3/1}(t)=
	\left\{\begin{array}{ll}
		\frac{(1-|t|)(1+(1+a_{01})|t|-a_{03}|t|^2)}{1+b_{01}|t|}, & 0 \leq |t|<1,\\
		\frac{(1-|t|)(2-|t|)(a_{10}+2a_{13}|t|)}{2(1+b_{11}|t|)}, & 1 \leq |t|<2,\\
		0,& |t|\geq 2.
	\end{array}\right.
\end{equation}

In the following, we examine the necessity for $S_{3/1}(t)$ to be a $C^1$ continuous function, which is equivalent to
\begin{equation*}
	\left\{
	\begin{aligned}
		&-a_{01}+b_{01}=a_{01}-b_{01},\\
		&\frac{-2-a_{01}+a_{03}}{1+b_{01}}=-\frac{2a_{13}+a_{10}}{2+2b_{11}},\\
		&\frac{4a_{13}+a_{10}}{2+4b_{11}}=0.
	\end{aligned}
	\right.
\end{equation*}
Solving these equations leads to
$a_{10}=-4a_{13}, b_{01}=a_{01}, a_{03}=(2+a_{13}+a_{01}+a_{13}a_{01}+2b_{11}+a_{01}b_{11})/(1+b_{11})$
under new premises $1+b_{11}\neq 0,(1+a_{01})(1+2b_{11})\neq 0$.
By substituting $a_{10}, b_{01}, a_{03}$ into equation (\ref{eqn:cubic/linear interpolation form}), we can obtain the $C^1$ continuous $S_{3/1}(t)$
in the form
\begin{equation*}
	S_{3/1}(t)=\left\{
	\begin{aligned}
		&\frac{(1-|t|) \left(1+(1+a_{01})|t|-a_{03} |t|^2\right)}{1+{a_{01}} |t|}, & 0\leq |t|<1,\\
		&\frac{a_{13}(-1+|t|)(2-|t|)^2}{1+b_{11}|t|}, & 1\leq |t|<2,\\
		& 0,~~&|t|\geq 2,
	\end{aligned}
	\right.
\end{equation*}
where $a_{03}= (2 + a_{13} + a_{01} + a_{13}a_{01} + 2 b_{11} + a_{01}b_{11})/(1 + b_{11}), a_{01}\neq 0, b_{12}\neq -1$ or  $-1/2$.
It can be easily verified that
if $a_{03}=2+a_{13}, a_{01}=b_{11}=0$, $S_{3/1}(t)$ degenerates to the cubic interpolation kernel $S_3(t)$ (\ref{eqn:cubic interpolation kernel}).


A good interpolation kernel should have the partition of unity property, which for $S_{3/1}(t)$ means satisfying
$$
S_{3/1}(t-2)+S_{3/1}(t-1)+S_{3/1}(t)+S_{3/1}(t+1)=1, \forall t\in [0,1],
$$
or equivalently,
	\begin{equation}
		\label{eqn: partition of unity property condition a1 e1 a2 e2}
		\begin{split}
			\frac{a_{13}(-1+t)t^2}{-1-2b_{11}+b_{11}t}+\frac{t(-2+a_{03}-a_{01}+(1-2a_{03}+a_{01})t+a_{03}t^2)}{-1-a_{01}+a_{01}t}\\
			+\frac{(1-t)(1+(1+a_{01})t-a_{03}t^2)}{1+a_{01}t}
			\;\;\;+\frac{a_{13}t(1-t)^2}{1+b_{11}+b_{11}t}=1,\forall t \in [0,1].
		\end{split}
	\end{equation}
After rewriting equation (\ref{eqn: partition of unity property condition a1 e1 a2 e2}) as a rational function that equals  zero  and requiring
that the numerator of the rational function to be zero, we can obtain the
following four independent equations
	\begin{equation}
		\label{eqn:several conditions unity partition}
		\left\{
		\begin{aligned}
			&\;\;\left(2 b_{11}+1\right) \left(a_{03} b_{11}-a_{13}a_{01} -a_{01} b_{11}+a_{03}-a_{01}-a_{13}-2 b_{11}-2\right)=0,\\
			&\;\;a_{13}a_{01}^2 +2 a_{13}a_{01}^2  b_{11}-4 a_{03} a_{01} b_{11}^2-a_{13}a_{01} -6 a_{03} a_{01} b_{11}-4 a_{13}a_{01} b_{11}\\
			& +a_{03} b_{11}^2+3 a_{03} b_{11}-a_{01} b_{11}^2-3 a_{01} b_{11}-2 a_{03} a_{01}+a_{03} -a_{01}-4 a_{13} b_{11}\\
			& -a_{13}-2 b_{11}^2-6 b_{11}-2=0,\\
			&\;\;a_{13}a_{01}^2 +8 a_{13} a_{01}^2 b_{11}+2 a_{03} a_{01} b_{11}^2-6 a_{03} a_{01} b_{11}-2a_{13} a_{01}  b_{11}-a_{03} b_{11}^2\\
			&+a_{01} b_{11}^2-2 a_{03} a_{01}-2 a_{13} b_{11}+2 b_{11}^2=0,\\
			& \;\;a_{01} b_{11} \left(a_{03} b_{11}+a_{13}a_{01}\right)=0.
		\end{aligned}
		\right.
	\end{equation}

Under the conditions $a_{03}= (2 + a_{13} + a_{01} + a_{13} a_{01} + 2 b_{11} + a_{01} b_{11})/(1 + b_{11})$, $a_{01}\neq-1, b_{11}\neq-1\; \text{or} \;-1/2$ and
the denominators of  $S_{3/1}(t)$ have no roots in corresponding intervals,  the equations (\ref{eqn:several conditions unity partition})
lead to six solutions:
		\begin{enumerate}[(i)]
			\item $a_{03}=2, a_{13}=a_{01}=0,$
			\item  $a_{03}=2, a_{13}=a_{01}=b_{11}=0,$
			\item  $a_{03}=-a_{01}, a_{13}=-2, b_{11}=0,$
			\item  $a_{03}=1, a_{13}=-1, a_{01}=b_{11}=0,$
			\item  $a_{03}=2+a_{13}, a_{01}=b_{11}=0,$
			\item  $a_{03}=1, a_{13}=-1-b_{11}, a_{01}=\frac{b_{11}}{1+b_{11}}.$
		\end{enumerate}

It is easy to check that the solution (i) to solution (iv) are all special cases of the solution (v), which is exactly  the well-known cubic interpolation kernel $S_3(t)$.
Thus, the four kernels corresponding to the solution (i) to solution (iv) are all special cases of the cubic interpolation kernel $S_3(t)$.
The solution (vi) leads to a new  cubic/linear rational kernel $S_{3/1}(t)$ with the form
\begin{equation}
	\label{eqn: 3/1 after partition with b Form}
	S_{3/1}(t)=\left\{
	\begin{aligned}
		&\frac{(1-|t|)(1+(1+a_{01})|t|-t^2)}{1+a_{01}|t|}, ~~~&& 0\leq |t|<1,\\
		&\frac{(1-|t|)(2-|t|)^2}{1-a_{01}+a_{01}|t|}, ~~~&& 1\leq |t|<2,\\
		&0, ~~~&& |t|\geq 2.
	\end{aligned}
	\right.
\end{equation}



During  the derivation of the interpolation kernel $S_{3/1}$ (\ref{eqn: 3/1 after partition with b Form}), we neglected  inequality requirements and the limitation that singularities should not exist in all denominators.
To ensure the rigor of the deduction,
we employ a direct verification approach instead of considering these neglected constraints.
Through computations, we demonstrate that  $S_{3/1}(t)$ (\ref{eqn: 3/1 after partition with b Form}) is  $C^1$ continuous, interpolates the sinc function at integer points,
and satisfies the partition of unity property directly.
To avoid singularities, only the condition $a_{01}\ge -1$ must be satisfied.

We will now present some common choices for the parameter $a_{01}$ in $S_{3/1}(t)$.
(a)
If $a_{01}=0$, $S_{3/1}(t)$ 
	degenerates to a spline interpolation kernel
which corresponds exactly to the case where
$a_{02}=-2$ in $S_3(t)$ (\ref{eqn:cubic interpolation kernel}).


(b)
If $a_{01}\rightarrow +\infty$, $S_{3/1}(t)$ becomes the linear interpolation kernel $S_1(t)$ (\ref{eqn: linear kernel}).

(c)
When $a_{01}=-1$, the common factor $1-t$ appears in both the numerators and denominators.
Removing the common factor $1-t$ in the numerators and denominators results in
\begin{equation}
	\label{eqn: S2 final form}
	S_{2}(t)=\left\{
	\begin{aligned}
		&1-|t|^2, ~~~&& 0\leq |t|<1,\\
		&(1-|t|)(2-|t|), ~~~&& 1\leq |t|<2,\\
		&0, ~~~&& |t|\geq 2.
	\end{aligned}
	\right.
\end{equation}
In this case, one can check that it has  approximation order 2.

\section{Rational quartic/linear interpolation kernel \(S_{4/1}(t)\)}
\label{section:Rational quartic/linear interpolation kernel $S_{4/1}(t)$}
In this section, we construct symmetric quartic/linear  rational interpolation kernels supported on $[-2,2]$.
The general form of this type of kernel is
\begin{equation*}
	S_{4/1}(t)=
	\left\{\begin{array}{ll}
		\frac{a_{00} +a_{01}|t|+a_{02}|t|^{2}+a_{03}|t|^{3}  +a_{04}|t|^{4}}{1 +b_{01}|t|}, & 0 \leq |t|<1,\\
		\frac{a_{10}+a_{11}|t|+a_{12}|t|^{2}  +a_{13}|t|^{3}  +a_{14}|t|^{4}}{1+b_{11}|t|}, & 1 \leq |t|<2,\\
		0, & |t|\geq 2.
	\end{array}\right.
\end{equation*}
Similarly,
$S_{4/1}(t)$ should interpolate the sinc function at $t=0, 1, 2$, i.e.,
$S_{4/1}(0)=\text{sinc}(0)=1,
S_{4/1}(1)=\lim_{t\to1^{-}}S_{4/1}(t)=\text{sinc}(1)=0,
S_{4/1}(2)=\lim_{t\to2^{-}}S_{4/1}(t)=\text{sinc}(2)=0$,
which are equivalent to
\begin{equation*}
	\left\{
	\begin{aligned}
		&a_{00}= 1,\\
		& \frac{a_{00}+a_{01}+a_{02}+a_{03}+a_{04}}{1+b_{01}} = 0,\\
		&\frac{a_{10}+a_{11}+a_{12}+a_{13}+a_{14}}{1+b_{11}} = 0,\\
		&\frac{a_{10}+2a_{11}+4a_{12}+8a_{13}+16a_{14}}{1+2b_{11}} = 0.
	\end{aligned}
	\right.
\end{equation*}
 This implies that  $b_{01}\neq -1, b_{11}\neq -1, b_{11}\neq-1/2, a_{00}=1$, and
\begin{equation*}
	\left\{
	\begin{aligned}
		&a_{04}=-1-a_{01}-a_{02}-a_{03},\\
		&a_{13}=-\frac{15}{8}a_{10}-\frac{7}{4}a_{11}-\frac{3}{2}a_{12},\\
		&a_{14}=\frac{7}{8}a_{10}+\frac{3}{4}a_{11}+\frac{1}{2}a_{12}.
	\end{aligned}
	\right.
\end{equation*}
It follows that  $S_{4/1}(t)$ is
\begin{equation}
	\label{eqn:4/1 after interpolation}
	\left\{\begin{array}{cc}
		\frac{1+a_{01}|t| +a_{02}|t|^{2} + a_{03}|t|^{3} +(-1-a_{01}-a_{02}-a_{03})|t|^{4}}{1+b_{01}|t|}, & 0\leq |t|<1,\\
		\frac{a_{10}+a_{11}|t| +a_{12}|t|^{2} +( -\frac{15}{8}a_{10}-\frac{7}{4}a_{11}-\frac{3}{2}a_{12})|t|^{3} +(\frac{7}{8}a_{10}+\frac{3}{4}a_{11}+\frac{1}{2}a_{12})|t|^{4}}{1+b_{11}|t|}, & 1 \leq |t|<2,\\
		0, & |t|\geq 2.
	\end{array}\right.
\end{equation}

It can be easily check that $C^1$ continuous conditions are equivalent to
\begin{equation*}
	\left\{
	\begin{aligned}
		&-a_{01}+b_{01}=a_{01}-b_{01},\\
		&\frac{4+3a_{01}+2a_{02}+a_{03}}{1+b_{01}}=\frac{17a_{10}+10a_{11}+4a_{12}}{8+8b_{11}},\\
		&\frac{11a_{10}+8a_{11}+4a_{12}}{2+4b_{11}}=0.
	\end{aligned}
	\right.
\end{equation*}
Solving these equations leads to  two solutions. The first solution is
\begin{equation}
	\label{eqn: derivative 0  4/1 solution}
	\left\{
	\begin{aligned}
		a_{03}&=-4-3a_{01}-2a_{02},\\
		a_{11}&=-3a_{10},\\
		a_{12}&=\frac{13}{4}a_{10},\\
		b_{01}&=a_{01}
	\end{aligned}
	\right.
\end{equation}
under the condition $1+a_{01}+3b_{11}+3a_{01}b_{11}+3a_{01}b_{11}+2b_{11}^2+2a_{01}b_{11}^2\neq 0$,
and the second one is
\begin{equation}
	\label{eqn: derivative  nonzero  4/1 solution}
	\left\{
	\begin{aligned}
		a_{12}&=\frac{1}{4}(-11a_{10}-8a_{11}),\\
		b_{01}&=a_{01},\\
		b_{11}&=\frac{-16-12a_{01}-8a_{02}-4a_{03}+3a_{01}a_{10}+a_{11}+a_{01}a_{11}}{4(4+3a_{01}+2a_{02}+a_{03})}
	\end{aligned}
	\right.
\end{equation}
under the conditions $4+3a_{01}+2a_{02}+a_{03}\neq 0$
and
	\begin{equation*}
		\begin{aligned}
			&\;\;-24a_{10}-42a_{01}a_{10}-18a_{01}^2a_{10}-12a_{02}a_{10}-12a_{01}a_{02}a_{10}-6a_{03}a_{10}\\
			&-6a_{01}a_{03}a_{10}+9a_{10}^2+18a_{01}a_{10}^2+9a_{01}^2a_{10}^2-8a_{11}-14a_{01}a_{11}-6a_{01}^2a_{11}\\
			&-4a_{02}a_{11}-4a_{01}a_{02}a_{11}-2a_{03}a_{11}-2a_{01}a_{03}a_{11}+6a_{10}a_{11}+12a_{01}a_{10}a_{11}\\
			&+6a_{01}^2a_{10}a_{11}+a_{11}^2+2a_{01}a_{11}^2+a_{01}^2a_{11}^2\neq 0.\\
		\end{aligned}
	\end{equation*}

\subsection{Quartic/linear interpolation kernel function with zero derivative at \(t=1\)}
\label{subsection:4/1 kernel function with zero derivative at t}

Substituting the first solution (\ref{eqn: derivative 0  4/1 solution}) into (\ref{eqn:4/1 after interpolation}) yields
	\begin{equation}
		\label{eqn: 1st C1}
		S_{4/1}(t)=
		\left\{
		\begin{aligned}
			&\frac{ (1-|t|)^2(1+(2+a_{01})|t|+(3+2a_{01}+a_{02})|t|^2  )}{1+a_{01}|t|},\! &0\leq |t| <1,\\
			&\frac{ a_{10}(1-|t|)^2(2-|t|)^2}{4(1+b_{11}|t|)},\! &1\leq |t| <2,\\
			&0, & |t|\geq 2.
		\end{aligned}
		\right.
	\end{equation}
If the corresponding numerators and denominators have no common factors, its derivative at $t=1$ is  0.

Next, we consider the partition of unity property. Similar to the derivation process of the interpolation kernel (\ref{eqn: 3/1 after partition with b Form}), we impose
the condition
$$
S_{4/1}(t-2)+S_{4/1}(t-1)+S_{4/1}(t)+S_{4/1}(t+1)=1, \forall t\in [0,1].
$$
This constraint leads to a complicated nonlinear  system of equations,
which we omit here for brevity.
We then solve these equations, and discard  solutions that lead to trivial polynomial spline interpolation kernels or rational interpolation kernels with singularities. Consequently, only three meaningful solutions are left.
These correspond to rational kernels with zero derivative at $t=1$, 
as illustrated below.

\subsubsection{The first quartic/linear interpolation kernel \(S^1_{4/1}(t)\) with zero derivative at \(t=1\)}

The first solution of the three meaningful solutions is $a_{1 0} =-4(3+a_{0 2})/(1+2 a_{0 1}), b_{11} =-a_{0 1}/(1+2 a_{0 1})$, we substitute it
into equation (\ref{eqn: 1st C1}) leads to the first quartic/linear interpolation kernel
\begin{equation*}
	S^1_{4/1}(t)=
	\left\{\begin{array}{cl}
		\frac{(1-|t|)^{2}\left(1+(2 + a_{0 1})|t|+(3+2 a_{0 1}+a_{0 2})|t|^{2} \right)}{1+a_{0 1}|t| }, & 0 \leq |t|<1,\\
		\frac{(2-|t|)^{2}(1-|t|)^{2}(3+a_{0 2})}{-1-2 a_{0 1}+a_{0 1}|t| }, & 1 \leq |t|<2,\\
		0, & |t|\geq 2\end{array}\right.
\end{equation*}
with zero derivative at $t=1$.

We did not consider any inequalities in the derivation process of $S^1_{4/1}$.
Similar to the derivation process of the interpolation kernel (\ref{eqn: 3/1 after partition with b Form}),
we also verify this kernel directly to ensure the rigor of our argument.
The computations demonstrate that the kernel interpolates the sinc function at integer points, exhibits $C^1$ continuity, and satisfies the partition of unity property. The remaining four $S_{4/1}$ type interpolation kernels undergo similar analysis, and they all satisfy these conditions directly.
Therefore, additional explanations are not necessary as these kernels should satisfy  similar constraints and can be analyzed similarly.

To avoid singularities in the two denominators, $a_{01}$ must satisfy the condition $a_{01}>-1$.
The common factor $1-t$ appears in all its numerators and denominators if $a_{01}=-1$, resulting in the degeneration of the kernel to the cubic interpolation kernel $S_3(t)$ if we remove the common factor $1-t$ from the numerators and denominators.
Analysis similar to that used in the derivation of $S_{3/1}(t)$ reveals that the kernel achieves the highest approximation order of 3 if $a_{01}=-1$ and $a_{02}=-5/2$, which is precisely when $S_3(t)$ attains the highest approximation order.

\subsubsection{The second quartic/linear interpolation kernel $S^2_{4/1}$ with zero derivative at $t=1$}

The second solution of the three meaningful solutions is $a_{1 0} =-4(3+a_{0 2})/(1- a_{0 1}), b_{11} =a_{0 1}/(1- a_{0 1})$.
By substituting it
into equation (\ref{eqn: 1st C1}), we derived the second quartic/linear  rational function with zero derivative at $t=1$ as following
\begin{equation*}
	S^2_{4/1}(t)=\left\{\begin{array}{cl}
		\frac{(1-|t|)^{2}\left(1+(2 + a_{0 1})|t|+(3+2 a_{0 1}+a_{0 2})|t|^{2} \right)}{1+a_{0 1}|t| }, & 0 \leq |t|<1,\\
		\frac{(2-|t|)^{2}(1-|t|)^{2}(3+a_{0 2})}{-1+ a_{0 1}-a_{0 1}|t| }, & 1 \leq |t|<2,\\
		0, & |t|\geq 2.\end{array}\right.
\end{equation*}

To avoid singularities in the two denominators,  it is easy to check that only $a_{01}$ should meet the condition $a_{01}\geq-1$.
If $a_{01}=-1$, we remove the common factor $1-t$ in numerators and
denominators,
then it degenerates into a  cubic interpolation kernel
	\begin{equation}
		\label{eqn:slightly different cubic interpolation kernel}
		\left\{
		\begin{aligned}
			&(1-|t|)(1+|t|+(1+a_{02})|t|^2), & 0\leq |t|<1,\\
			&-(3+a_{02})(2-|t|)(1-|t|)^2\!,& 1\leq |t|<2,\\
			&0, &  |t|\geq 2.
		\end{aligned}
		\right.
	\end{equation}
This kernel has same expression as $S_3(t)$ on the interval $(-1,1)$, but differs slightly on the intervals $(-2,-1)$ and $(1,2)$.

If $a_{01}=-1$ and  $a_{02}=-4$, $S^2_{4/1}(t)$ degenerates into
\begin{equation*}
	\left\{\begin{array}{cl}
		(1-|t|)\left(1+|t|-3|t|^{2} \right), & 0 \leq |t|<1,\\
		(2-|t|)(1-|t|)^{2}, & 1 \leq |t|<2,\\
		0, & |t|\geq 2.\end{array}\right.
\end{equation*}
In this case,  we can easily verify by definition that it has the highest approximation order of 2.

\subsubsection{The third quartic/linear interpolation kernel $S^3_{4/1}$with zero derivative at  $t=1$}

The last solution of the three meaningful solutions is $a_{1 0} =-8(3+a_{0 2})/3, b_{11} =-1/3$.
Substituting
it 
into equation (\ref{eqn: 1st C1}) leads to the third quartic/linear interpolation kernel with zero derivative at $t=1$ as following
\begin{equation*}
	S_{4/1}^3(t)=\left\{\begin{array}{cl}
		\frac{(1-|t|)^{2}\left(2+3|t|+(2 a_{0 1}+4)|t|^{2} \right)}{2-|t| }, & 0 \leq |t|<1,\\
		\frac{(2-|t|)^{2}(1-|t|)^{2}(6+2a_{0 2})}{-3+|t| }, & 1 \leq |t|<2,\\
		0, & |t|\geq 2.\end{array}\right.
\end{equation*}
Though it has two free parameters, we can check by the definition (\ref{eqn:approximation order definition}) that  it cannot have higher approximation order than 1.

\subsection{Quartic/linear interpolation kernel with nonzero derivative at $t=1$}
By substituting the second solution (\ref{eqn: derivative  nonzero  4/1 solution}) of $C^1$ constraints into (\ref{eqn:4/1 after interpolation}), we obtain
	\begin{equation*}
		S_{4/1}(t)=\left\{
		\begin{aligned}
			&\frac{ (1-|t|)(1+(1+a_{01})|t|+(1+a_{01}+a_{02})|t|^2 +(1+a_{01}+a_{02}+a_{03})|t|^3 )}{1+a_{01}|t|},\\
			&\qquad\qquad\qquad\qquad\qquad \qquad\qquad\qquad\qquad\qquad\qquad 0\leq |t| <1,\\
			&\frac{ (1-|t|)(2-|t|)^2(a_{10}+(2a_{10}+a_{11})|t|)}{4(1+b_{11}|t|)}, \qquad\qquad\qquad~ 1\leq |t| <2,\\
			&0, \qquad\qquad\qquad\qquad\qquad\qquad\qquad\qquad\qquad\qquad\qquad|t|\geq 2,
		\end{aligned}
		\right.
	\end{equation*}
where $	b_{11}=\frac{-16-12a_{01}-8a_{02}-4a_{03}+3a_{01}a_{10}+a_{11}+a_{01}a_{11}}{4(4+3a_{01}+2a_{02}+a_{03})}$.
The value of its derivative at $t=1$ is $-(4+3a_{01}+2a_{02}+a_{03})/(1+a_{01})$.

Similar to the derivation process of the interpolation kernels in subsection \ref{subsection:4/1 kernel function with zero derivative at t}, we also consider the partition of unity property by definition.
This also leads to a complicated nonlinear system of equations, which we will not present here.
We then solve these equations and discard the solutions that result in trivial polynomial spline kernels or rational interpolation kernels with singularities. Ultimately, only two meaningful solutions remain, which are introduced below.

\subsubsection{The first quartic/linear interpolation kernel $S_{4/1}^4$ with nonzero derivative at $t=1$}

The first meaningful solution is $
a_{1 0} = \frac{4\left(-5+a_{0 1}+3 a_{0 1}^{2}-3 a_{0 2}+3 a_{0 1} a_{0 2}-2 a_{0 3}+a_{0 1} a_{0 3}\right)}{-1+a_{0 1}^{2}}$ and
$a_{11}=\frac{4\left(-11+6 a_{0 1}+9 a_{0 1}^{2}-7 a_{0 2}+9 a_{0 1} a_{0 2}-5 a_{0 3}+3 a_{0 1} a_{0 3}\right)}
{1-a_{0 1}^{2}}.
$ Substituting it
into equation (\ref{eqn: 1st C1}), we obtain the first interpolation kernel a nonzero derivative at $t=1$, which takes the form
	\begin{equation}
		\label{eqn:first 4/1 nonzero}
		S_{4/1}^4(t)=\left\{
		\begin{aligned}
			&\frac{ (1-|t|)\left(1+(1+a_{01})|t|+(1+a_{01}+a_{02})|t|^2  +(1+a_{01}+a_{02}+a_{03})  |t|^3  \right) }{1+a_{01}|t|},\\
			&\qquad\qquad\qquad\qquad\qquad\qquad\qquad\qquad\qquad\qquad~~0\leq |t| <1,\\
			&\frac{ (1-|t|)(2-|t|)^2(A+B|t|)}{(1+a_{01})(1-a_{01}+a_{01}|t|)}, \qquad\qquad\qquad\qquad\qquad1\leq |t| <2,\\
			&0,\qquad\qquad\qquad\qquad\qquad\qquad\qquad\qquad\qquad\qquad |t|\geq 2,
		\end{aligned}
		\right.
	\end{equation}
where $A=5-a_{01}-3a_{01}^2+3a_{02}-3a_{01}a_{02}+2a_{03}-a_{01}a_{03}$ 
and $B=-1+4a_{01}+3a_{01}^2-a_{02}+3a_{01}a_{02}-a_{03}+a_{01}a_{03}$.

To avoid singularities in the denominators, $a_{01}$ should satisfy the condition $a_{01}>-1$.
By selecting $a_{02}=-2-a_{01}$ and $a_{03}=1$, the interpolation kernel (\ref{eqn:first 4/1 nonzero}) turns into
the rational  interpolation kernel $S_{3/1}(t)$.
Choosing $a_{03}=-1-a_{02}+a_{01}a_{02}$ results in the cubic interpolation kernel $S_3(t)$.
When $a_{01}=0$, the denominators disappear, yielding a quartic polynomial
	\begin{equation}
		\label{eqn:quartic polynomial kernel}
		S_{4}(t)=\left\{
		\begin{aligned}
			& (1-|t|)\left(1+|t|+(1+a_{02})|t|^2  +(1+a_{02}+a_{03})  |t|^3  \right),& 0\leq |t| <1,\\
			& (1-|t|)(2-|t|)^2(5+3a_{02}+2a_{03}-(1+a_{02}+a_{03})|t|), &1\leq |t| <2,\\
			&0, & |t|\geq 2.
		\end{aligned}
		\right.
	\end{equation}
It is easy to verify that if $a_{03}=-1-a_{02}$, the quartic polynomial interpolation kernel $S_{4}(t)$ degenerates into the cubic interpolation kernel $S_3(t)$.

If $a_{03}=(-7-4a_{02})/2$, then $S_4(t)$ (\ref{eqn:quartic polynomial kernel}) has approximation order 2 with the free parameter $a_{02}$.
When $a_{02}=-5/2, a_{03}=3/2$, $S_{4}(t)$  has the highest approximation order 3,
and  it is exactly the cubic interpolation kernel $S_3(t)$ with $a_{02}= -5/2$.
One advantage of the interpolation kernel $S_4(t)$ (\ref{eqn:quartic polynomial kernel}) is that it can have approximation order 2 with a free parameter,
 while  $S_3(t)$ still has no free parameter since when $S_3(t)$  is required to have   approximation order 2, the parameter $a_{02}$ is also taken to be $-5/2$.

If $a_{02}=-6(3+a_{01})/(6+a_{01}), a_{03}=(-36-18a_{01}-17a_{02})/6$, then it is $C^2$ continuous.
In addition, if $a_{02}=-3(94+2a_{01}-4a_{01}^2+3a_{01}^3)/136$ and $a_{01}$ is the real root of the polynomial $-12+10x+26x^2+14x^3+3x^4$,
then it is $C^3$ continuous.

\subsubsection{The second quartic/linear interpolation kernel $S_{4/1}^4$ with nonzero derivative at $t=1$}
Substituting the second meaningful solution
$a_{1 0} = \frac{4(5+6 a_{0 1}+3 a_{0 2}+2 a_{0 3})}{1+2 a_{0 1}}$,
$a_{11} =-\frac{4(11+15 a_{0 1}+7 a_{0 2}+5 a_{0 3})}{1+2 a_{01}}$
into equation (\ref{eqn: 1st C1})  yields the second quartic/linear interpolation kernel  with nonzero derivative at $t=1$
	\begin{equation}
		\label{eqn:second 4/1 nonzero}
		S^5_{4/1}(t)=\left\{\begin{array}{cc}
			\frac{(1-|t|)\left(1+(1+a_{0 1})|t|+(1+a_{0 1}+a_{0 2})|t|^2+(1+a_{0 1}+a_{0 2}+a_{0 3})|t|^3\right)}{1+a_{0 1}|t|},& 0 \leq |t|<1,\\
			\frac{(1-|t|)(2-|t|)^2\left(5+6 a_{0 1}+3 a_{0 2}+2 a_{0 3}-(1+3 a_{0 1}+a_{0 2}+a_{0 3})|t|\right)}{1+2 a_{0 1}- a_{0 1}|t|},& 1 \leq |t|<2,\\
			0, & |t|\geq 2.
		\end{array}\right.
	\end{equation}

Similar to the interpolation kernel $S_{4/1}^4(t)$ (\ref{eqn:first 4/1 nonzero}),
no singularities occur in the two denominators of the interpolation kernel (\ref{eqn:second 4/1 nonzero}) if $a_{01}$ satisfies the condition $a_{01}>-1$.

After setting $a_{01}$ to zero, the interpolation kernel reduces to the quartic function $S_4(t)$ (\ref{eqn:quartic polynomial kernel}).
By selecting $a_{03}$ to be $-1-a_{02}+a_{01}a_{02}$, it becomes the cubic interpolation kernel $S_3(t)$ (\ref{eqn:cubic interpolation kernel}).

One can verify that if its approximation order is greater than 1, the interpolation kernel cannot remain in rational form.
In fact, it degenerates into the quartic polynomial $S_4(t)$ (\ref{eqn:quartic polynomial kernel}) if its approximation order exceeds 1.
Thus, the parameter selection method for approximation order 2 and 3 is the same as that for $S_4(t)$, except that $a_{01}=0$ is added.
In addition, if $a_{02}=-3, a_{03}=\frac{15-4a_{01}-9a_{01}^2}{2(3+2a_{01})}$, then it is $C^2$ continuous.

\begin{table*}[t]
	\centering
	\caption{Properties of all interpolation kernels.}
	\resizebox{0.95\columnwidth}{!}
	{
		\begin{tabular}{lcccccccccc}
			\toprule
			Kernels  & Support & Parameters & Continuity &Derivative ($t=1$) & $\int_{-\infty}^{\infty}S(\cdot )\; d\cdot$ & Approximation order & Special cases\\
			\midrule
			$S_0$ & $[-1/2,1/2]$&0&$C^{-1}$ &/& 1&1&/\\
			$S_1$ & $[-1,1]$&0&$C^0$&/&1&2&/\\
			$S_3$ & $[-2,2]$&1&$C^1$&$-3-a_{02}$ &1&3&/\\
			$S_{3/1}$&  $[-2,2]$&1&$C^1$&$-1$ &1&2&$S_1, S_2$\\
			$S^1_{4/1}$&$[-2,2]$&2&$C^1$&0 &1&3&$S_3$\\
			$S^2_{4/1}$&$[-2,2]$&2&$C^1$&0&1&2&/\\
			$S^3_{4/1}$&$[-2,2]$&1&$C^1$& 0&1&1&/\\
			$S^4_{4/1}$&$[-2,2]$&3&$C^3$&$-\frac{4+3a_{01}+2a_{02}+a_{03}}{1+a_{01}}$&1&3&$S_{3/1}$,$S_4$,$S_3$\\
			$S^5_{4/1}$&$[-2,2]$&3&$C^2$&$-\frac{4+3a_{01}+2a_{02}+a_{03}}{1+a_{01}}$&1&3&$S_4$,$S_3$\\
			\bottomrule
		\end{tabular}
	}
	\label{tab:Properties of all  kernels}
\end{table*}
\subsection{Summary of quartic/linear  interpolation kernels}
Table \ref{tab:Properties of all kernels} presents a comprehensive summary of the main properties of these interpolation kernels. The table includes the following columns:
\begin{itemize}
	\item \textbf{Parameters}: All five quartic/linear interpolation kernels have one or more free parameters, which may be adjusted to achieve desired properties in applications.
    \item \textbf{Continuity}: All five quartic/linear interpolation kernels are $C^1$ continuous. For certain parameter values, $S^4_{4/1}$ is $C^3$ continuous, and  $S^5_{4/1}$ can achieve $C^2$ continuity.
	\item \textbf{Derivative ($t=1$)}: $S^1_{4/1}$, $S^2_{4/1}$, and $S^3_{4/1}$ have zero derivative at $t=1$, while $S^4_{4/1}$ and $S^5_{4/1}$ have nonzero derivative $-\frac{4+3a_{01}+2a_{02}+a_{03}}{1+a_{01}}$ at $t=1$.
	\item \textbf{$\int_{-\infty}^{\infty}S(\cdot )\; d\cdot$}: All five quartic/linear interpolation kernels satisfy the partition of unity property directly, resulting in an integral value of 1 over their supports.
	\item \textbf{Approximation order}: $S^3_{4/1}$'s highest approximation order is 1, $S^2_{4/1}$'s highest approximation order is 2, and $S^1_{4/1}$, $S^4_{4/1}$, $S^5_{4/1}$'s highest approximation order is 3. However, all five interpolation kernels can only achieve approximation order 1 if they maintain their rational function forms. When the highest possible approximation order of these interpolation kernels is 3, they all degenerate into the cubic interpolation kernel $S_3(t)$ (\ref{eqn:cubic interpolation kernel}) with $a_{02}=-5/2$.
	\item \textbf{Special cases}: The interpolation kernel $S^1_{4/1}$, $S^4_{4/1}$, $S^5_{4/1}$ encompass $S_3$ as a particular instance. Additionally, $S^4_{4/1}$ and $S^5_{4/1}$ include $S_4$ as a special case, and $S_{3/1}$ is also a special case of $S^4_{4/1}$.

\end{itemize}

\section{Numerical experiments}
\label{sec:experiments}

\begin{table*}[t]
	\centering
	\caption{PSNR made by various interpolation kernels.}
	\resizebox{0.95\columnwidth}{!}
	{
		\begin{tabular}{lcccccccccc}
			\toprule
			Images  & Baboon & Barbara & Bike & Boats& Cameraman& Hat& Parrots& Pentagon& Peppers& Straw\\
			\midrule
			Nearest &  18.4118      &22.0045&19.2654&   22.4668   &   21.5854 & 25.6885 &   23.6639& 21.8686& 21.8422 & 18.5147\\
			Linear &18.5556&22.4552&19.8013&23.4429&22.1863& 26.3345& 24.3380&22.2607&22.9738 &18.6608 \\
			$S_3$ (Best $a_{02}$)&\makecell[c]{18.7618\\ (-1.805)} & \makecell[c]{22.7464\\(-2.175)} & \makecell[c]{20.4352 \\ (-1.965)} & \makecell[c]{24.3737 \\ (-1.960)} & \makecell[c]{22.7384\\ (-2.005)} & \makecell[c]{26.8300 \\ (-2.120)} & \makecell[c]{25.1072 \\ (-1.990)} & \makecell[c]{22.9283 \\ (-1.850)} & \makecell[c]{23.8671 \\ (-2.035)} & \makecell[c]{19.1962 \\ (-1.810)}\\
			$Hu\; and\; Tan's  \; S_{2/2}$ &18.7637 & 22.7582 & 20.4516 & 24.3964 & 22.7480 & 26.8452 & 25.1410 & 22.9315 & 23.9057 & 19.1877
			\\
			$S^4_{4/1}$(30,20,-121.5512)	& 18.7876	 &  22.7917&   20.5335 &  24.5061  & 22.8220 &  26.9177&25.2749 &   23.0157&  23.9955&  19.2776\\
			$S^4_{4/1}$(80,100,-444.7992) &\textbf{18.8034} & \textbf{22.7941} & \textbf{20.5662} & \textbf{24.5456} & \textbf{22.8449} & \textbf{26.9365} & \textbf{25.3517} & \textbf{23.0619} & \textbf{24.0055} & \textbf{19.3343}\\
			$S^5_{4/1}$(30,10,-90.1572) &18.7761&   22.7784 &20.4902&24.4480 &22.7740&26.8778 & 25.2105 &22.9739 & 23.9391& 19.2317\\
			$S^5_{4/1}$(50,10,-129.3052)&18.7725& 22.7706 & 20.4835& 24.4428& 22.7740&26.8708 &25.1922&22.9643&23.9377 & 19.2216\\
			\bottomrule
		\end{tabular}
	}
	\label{tab:best PSNR comparisions table For all kernels}
\end{table*}

\begin{table*}[t]
	\centering
	\caption{SSIM made by four interpolation methods.}
	\resizebox{1\columnwidth}{!}
	{
		\begin{tabular}{lcccccccccc}
					\toprule
					Images  & Baboon & Barbara & Bike & Boats& Cameraman& Hat& Parrots& Pentagon& Peppers& Straw\\
					\midrule
					Nearest &  0.2156      &0.4123&0.6283&  0.3824   &   0.2736 & 0.9319 &   0.8278& 0.3549& 0.4409 & 0.3034\\
					Linear &0.2029&0.4644&0.6478& 0.4128& 0.2729&0.9368 &0.8521 &0.3509 & 0.4964 & 0.2697\\
					$S_3$ (Best $a_{02}$)&\makecell[c]{0.2615\\ (-0.710)} & \makecell[c]{0.5006\\ (-2.170)} & \makecell[c]{0.6880 \\ (-1.920)} & \makecell[c]{0.4653 \\ (-1.930)} & \makecell[c]{0.3101 \\ (-1.930)} & \makecell[c]{0.9402 \\ (-2.430)} & \makecell[c]{0.8659 \\ (-2.255)} & \makecell[c]{0.4443 \\ (-1.270)} & \makecell[c]{0.5400 \\ (-2.145)} & \makecell[c]{0.3836\\ (-0.800)}\\
					$Hu \;and\; Tan's  \; S_{2/2}$ &0.2490 & 0.5045 & 0.6898 & 0.4680 & 0.3109 & 0.9393 & 0.8659 & 0.4346 & 0.5426 & 0.3587\\
					$S^4_{4/1}$(30,20,-121.5512)&0.2580 &  0.5116    &0.6966&   0.4776 &   0.3198  &  0.9402 &  0.8688 &  0.4487  &  0.5500 &   0.3761\\
					$S^4_{4/1}$(80,100,-444.7992) &	0.2686 & 0.5184 & 0.7022 & 0.4862 & 0.3270 & 0.9402 & 0.8704 & 0.4641 & 0.5551 & 0.3956\\
					$S^5_{4/1}$(30,10,-90.1572) &0.2553&0.5108 &0.6943&0.4745& 0.3154& 0.9396&0.8673&0.4445&0.5473&0.3703\\
					$S^5_{4/1}$(50,10,-129.3052)&0.2530 &0.5081& 0.6928&  0.4724&0.3142& 0.9396&0.8670& 0.4410& 0.5459& 0.3662\\
					\bottomrule
				\end{tabular}
	}
	\label{tab:best SSIM comparisions table For all kernels}
\end{table*}

\begin{table*}[t]
	\centering
	\caption{FSIM made by various interpolation kernels.}
	\resizebox{1\columnwidth}{!}
	{
		\begin{tabular}{lcccccccccc}
			\toprule
			Images  & Baboon & Barbara & Bike & Boats& Cameraman& Hat& Parrots& Pentagon& Peppers& Straw\\
			\midrule
			Nearest &  0.5836&  0.6922 & 0.6346 & 0.6908& 0.6910 & 0.7473 & 0.7920& 0.6244& 0.7207& 0.5923\\
			Linear &0.5859   & 0.7489  &  0.7022 &   0.7800   & 0.7371  &  0.7878    &0.8678 &   0.6789 &   0.8164 &   0.5798\\
			$S_3$ (Best $a_{02}$)&\makecell[c]{\textbf{0.7178}\\ (0.735)} & \makecell[c]{0.7816\\ (-1.855)} & \makecell[c]{0.7371 \\ (-1.905)} & \makecell[c]{0.8007 \\ (-2.185)} & \makecell[c]{0.7476 \\ (-2.355)} & \makecell[c]{0.8014 \\ (-2.275)} & \makecell[c]{0.8816 \\ (-2.295)} & \makecell[c]{0.7329 \\ (-1.325)} & \makecell[c]{0.8306 \\ (-2.325)} & \makecell[c]{0.7213 \\ (-0.170)}\\
			$Hu \; and\; Tan's  \; S_{2/2}$ & 0.6514 & 0.7846 & 0.7381 & 0.8007 & 0.7456 & 0.8013 & 0.8820 & 0.7267 & 0.8298 & 0.6629\\
			$S^4_{4/1}$(30,20,-121.5512)&0.6622 & 0.7900  & 0.7443 &0.8045 & 0.7489 &0.8050 & 0.8855 & 0.7351 & 0.8328 &0.6779\\
			$S^4_{4/1}$(80,100,-444.7992)&0.6756 & \textbf{0.7962} & \textbf{0.7494} & \textbf{0.8063} & \textbf{0.7502} & \textbf{0.8075} & \textbf{0.8881} & \textbf{0.7433} & \textbf{0.8340} & 0.6944\\
			$S^5_{4/1}$(30,10,-90.1572) & 0.6584&0.7892& 0.7417& 0.8029&0.7470& 0.8037& 0.8846&0.7318 &0.8322& 0.6718\\
			$S^5_{4/1}$(50,10,-129.3052)&0.6562 &0.7873 &0.7407&0.8026 & 0.7468& 0.8028& 0.8837& 0.7304&0.8313& 0.6693\\
			\bottomrule
		\end{tabular}
	}
	\label{tab:best FSIM comparisions table For all kernels}
\end{table*}

\begin{figure*}[h!]
	\centering
	\subfigure[Baboon]{
		\label{Fig.sub.baboon}
		\includegraphics[width=.9in]{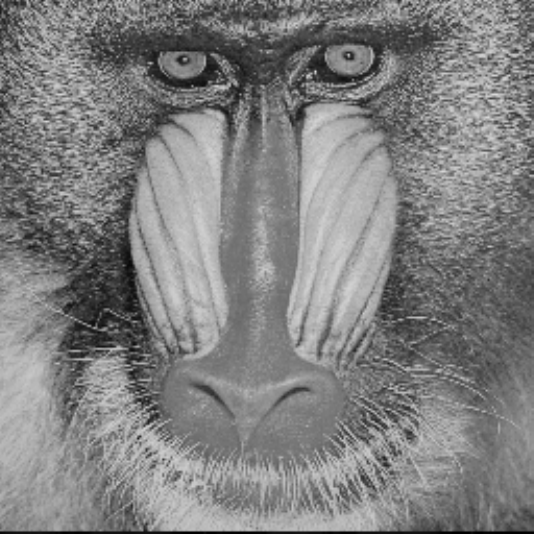}
	}
	\subfigure[Barbara]{
		\label{Fig.sub.barbara}
		\includegraphics[width=.9in]{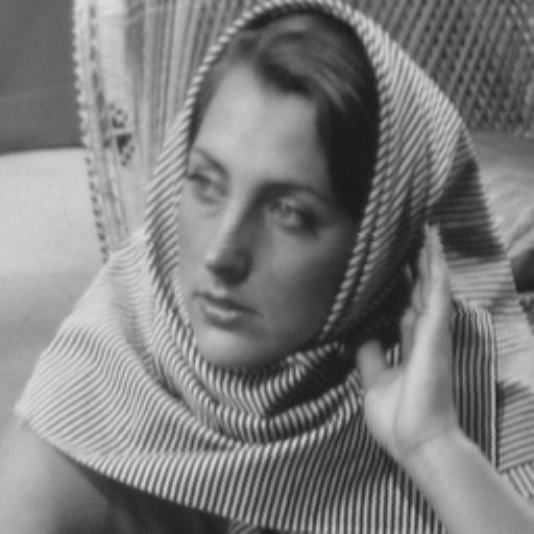}
	}
	\subfigure[Bike]{
		\label{Fig.sub.bike}
		\includegraphics[width=.9in]{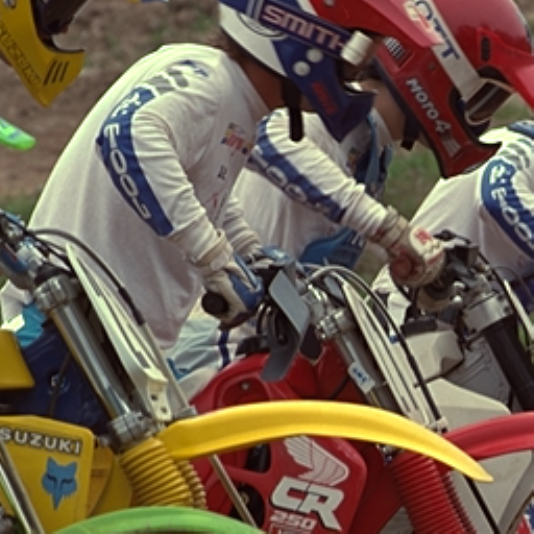}
	}
	\subfigure[Boats]{
		\label{Fig.sub.boats}
		\includegraphics[width=.9in]{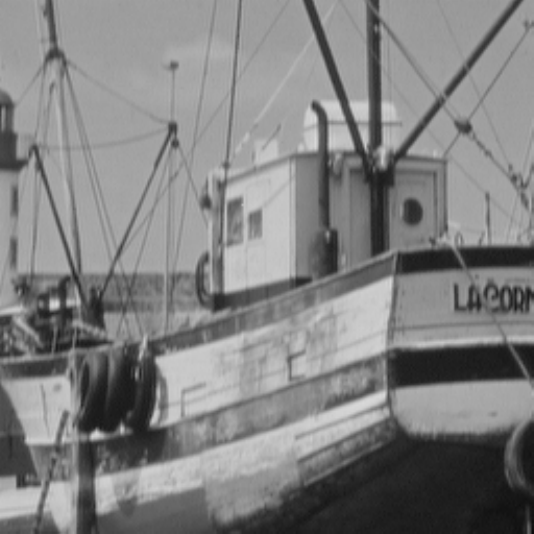}
	}
	\subfigure[Cameraman]{
		\label{Fig.sub.cameraman}
		\includegraphics[width=.9in]{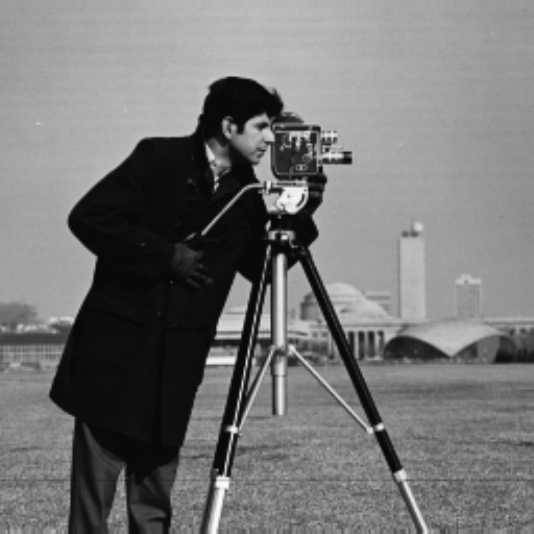}
	}
	\subfigure[Hat]{
		\label{Fig.sub.hat}
		\includegraphics[width=.9in]{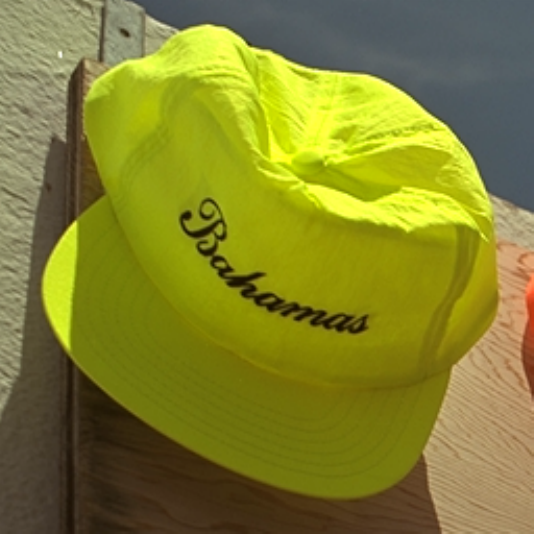}
	}
	\subfigure[Parrots]{
		\label{Fig.sub.Parrots}
		\includegraphics[width=.9in]{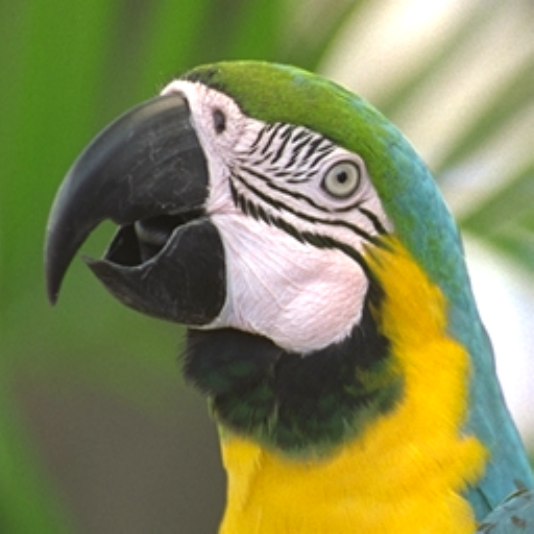}
	}
	\subfigure[Pentagon]{
		\label{Fig.sub.pentagon}
		\includegraphics[width=.9in]{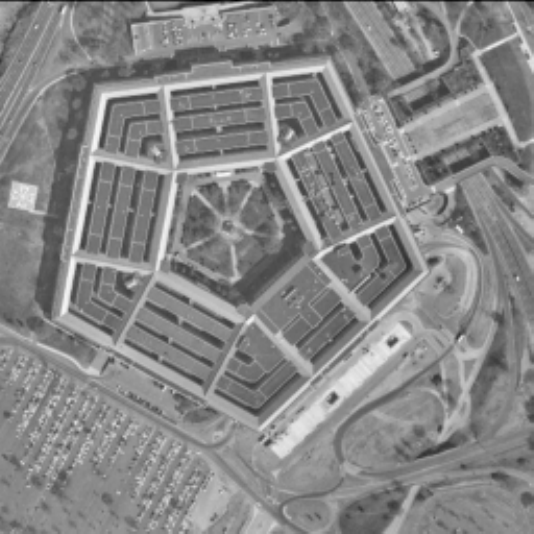}
	}
	\subfigure[Peppers]{
		\label{Fig.sub.peppers}
		\includegraphics[width=.9in]{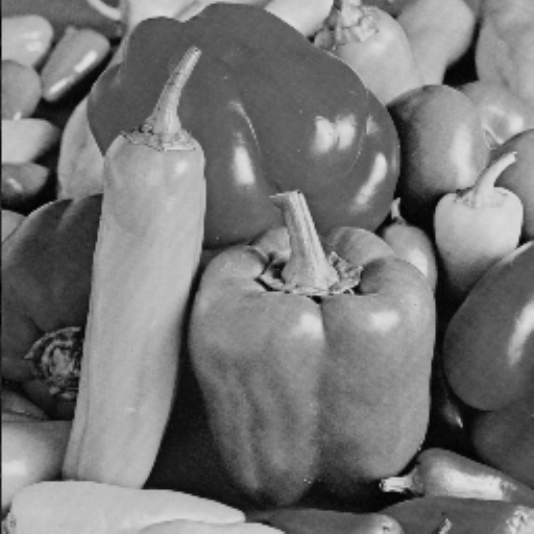}
	}
	\subfigure[Straw]{
		\label{Fig.sub.straw}
		\includegraphics[width=0.9in]{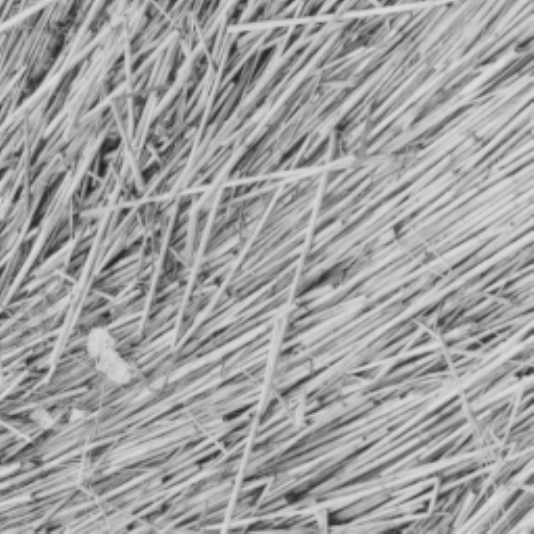}
	}
	\caption{
		Ten test images.
	}
	\label{Fig.ten test images}
\end{figure*}

Each of the interpolation kernels has one or more parameters, making it challenging to determine the optimal parameters.
In the following, we provide suggestions for selecting the parameters in $S^4_{4/1}(t)$ and $S^5_{4/1}(t)$ based on our experiments with image magnification. First, we selected ten 256$\times$256 images with dynamic ranges from $0$ to $255$, as shown in Figure \ref{Fig.ten test images}.

The ability of the interpolation kernels will be evaluated through image magnification, aiming to minimize the difference between the magnified image and the corresponding ideally magnified one, which unfortunately does not exist. To conduct the evaluation, all images are first reduced by a factor of 4$\times$4 using bicubic interpolation, and then magnified by a factor of 4.
	We utilized the built-in MATLAB function \textit{imresize} for  image magnification.
For image reduction, we used the bicubic parameter in \textit{imresize} function. In the case of image magnification, we replaced the interpolation kernel function used by \textit{imresize} with  custom interpolation kernels, specifically $S_3$, $S_{2/2}$, $S^4_{4/1}$ and $S^5_{4/1}$.
The evaluation indices used are peak signal-to-noise ratio (PSNR), structural similarity index (SSIM) \cite{wang2004imageSSIM}  and feature similarity index (FSIM) \cite{zhang2011FSIM}.
The nearest interpolation neighbor $S_0(t)$ (\ref{eqn:nearest neighbor kernel}), 
linear interpolation kernel $S_1(t)$ (\ref{eqn: linear kernel}), and cubic interpolation kernel $S_3(t)$ (\ref{eqn:cubic interpolation kernel}) were tested. For the cubic interpolation kernel $S_3(t)$, the range of $a_{02}$ was restricted to 
$[-7,1]$, and a step size of $0.005$ was selected. The highest PSNR values for the three types of interpolation kernels were recorded for each test image.
Table  \ref{tab:best PSNR comparisions table For all kernels} displays the optimal PSNR values for the ten test images in its first three rows.
The cubic interpolation kernel $S_3(t)$ exhibits varying the best parameters across the images.

\begin{figure*}[h!]
	\centering
	\includegraphics[width=2.5in]{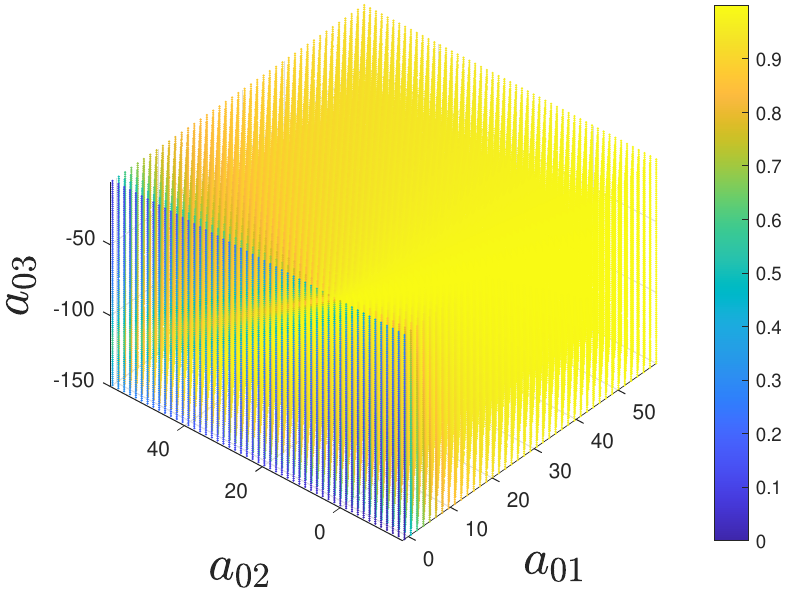}
	\includegraphics[width=2.5in]{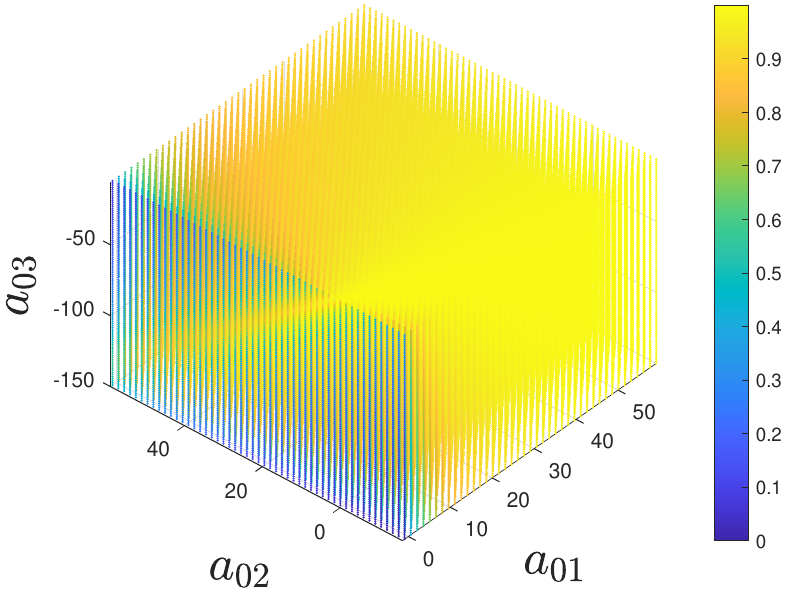}
	\caption{
		Normalized PSNR values for the ten test images (Left $S^4_{4/1}$, Right $S^5_{4/1}$).
	}
	\label{Fig: Normalized PSNR values for the ten test images}
\end{figure*}

To determine suitable ranges of free parameters,
we constructed an exhaustive search by iterating over the parameter ranges as follows: $a_{01}$ in $[-0.9, 70]$, $a_{02}$ in $[-16, 45]$, and $a_{03}$ in $[-150, -5]$, with a step size of 1.
Next, we applied a fixed set of parameters $(a_{01}, a_{02}, a_{03})$ to image magnification, and used the value of PSNR as the criterion for evaluating the quality of these parameters.
We enlarged the $10$ reduced test images by $S^4_{4/1}$ or $S^5_{4/1}$ for selected triplets $(a_{01}, a_{02}, a_{03})$ and computed corresponding PSNR values.
Due to the substantial discrepancies in PSNR values among diverse images, directly aggregating the PSNR values of ten images corresponding to a particular triplet $(a_{01}, a_{02}, a_{03})$ and employing this cumulative PSNR to assess the quality of a horizontal triplet is not a valid approach. Instead, we initially normalized all PSNR values of each image to the range $[0, 1]$ via unity-based normalization. We then computed the average of the normalized PSNR values corresponding to each parameter set. This methodology enables us to evaluate the comprehensive applicability of a specific parameter set across different images.
Figure \ref{Fig: Normalized PSNR values for the ten test images} illustrates the final averaged normalized PSNR values. The normalized PSNR value distributions of $S^4_{4/1}$  (Left part of Figure \ref{Fig: Normalized PSNR values for the ten test images}) and  $S^5_{4/1}$ (Right part of Figure \ref{Fig: Normalized PSNR values for the ten test images}) are highly similar. Both indicate that for small $a_{01}$ values, the high PSNR value zone is more concentrated. Conversely, for relatively large $a_{01}$ values, the high value region is less concentrated, and the overall value of the region is comparatively large.

\begin{figure*}[h!]
	\centering
	\subfigure[]{
		\label{Fig.sub.FirstN4D1 Better than cubic}
		\includegraphics[width=0.45\textwidth]{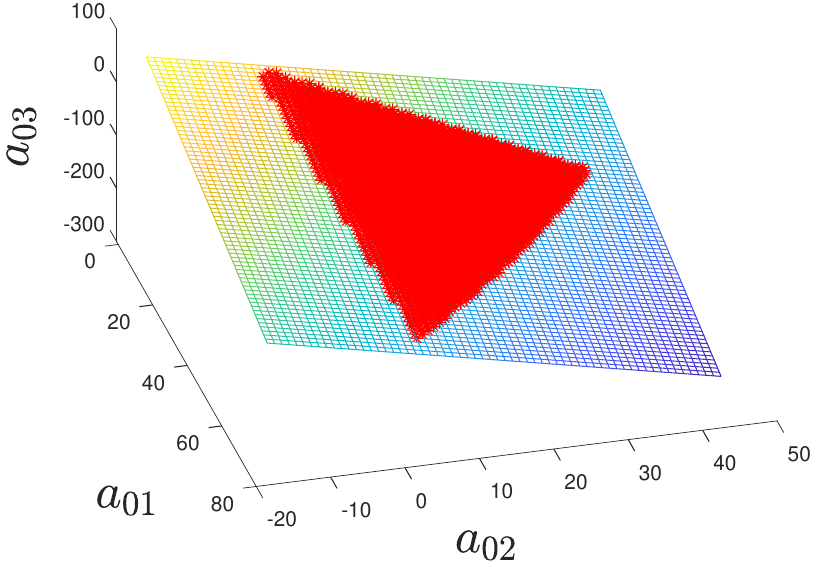}
	}
	\subfigure[]{
		\label{Fig.sub.FirstN4D1 Better than cubic a03 thickness angle}
		\includegraphics[width=0.45 \textwidth]{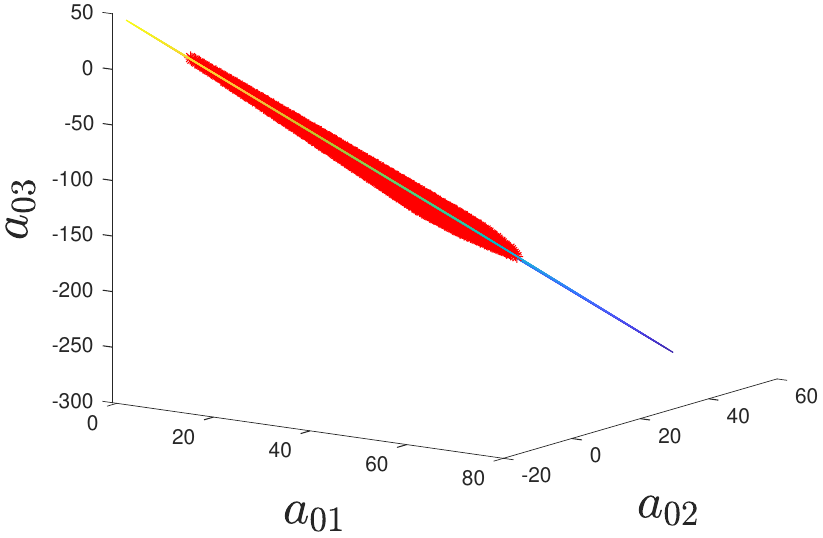}
	}
	\subfigure[]{
		\label{Fig.sub.SecondtN4D1 Better than cubic}
		\includegraphics[width=0.45 \textwidth]{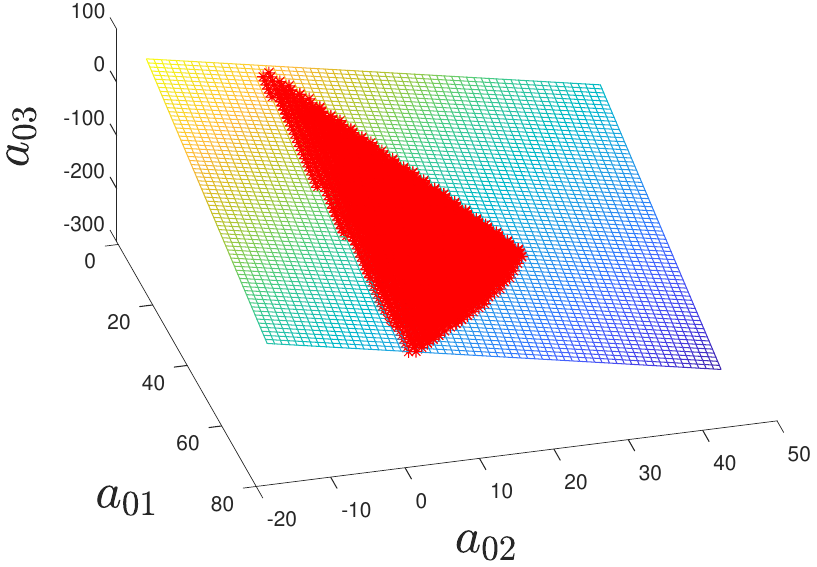}
	}
	\subfigure[]{
		\label{Fig.sub.SecondN4D1 Better than cubic a03 thickness angle}
		\includegraphics[width=0.45 \textwidth]{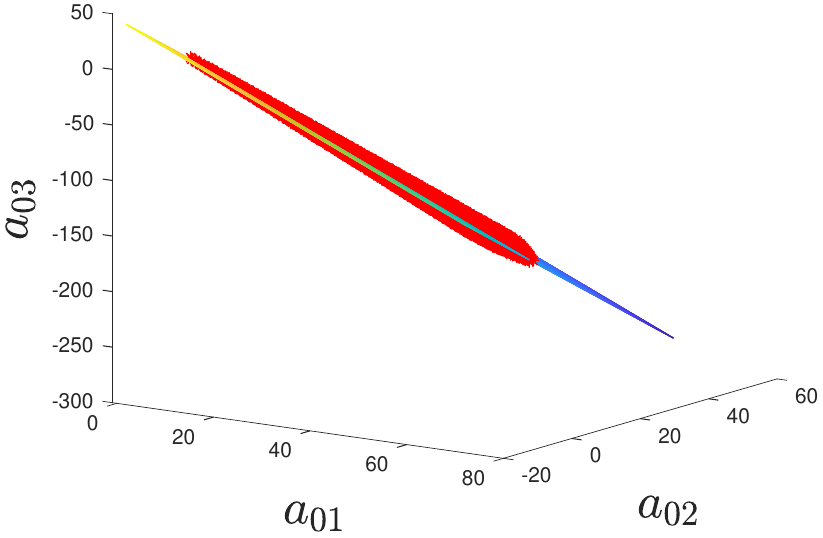}
	}
	\caption{Parameter distributions of $S_{4/1}^4$(First row) and $S_{4/1}^5$ (Second row) corresponding better PSNR values than those of  three kinds of interpolation kernels.}
	\label{Fig.First N4D1 better than cubic}
\end{figure*}

We  directly compare $S^4_{4/1}$ and $S^5_{4/1}$   with the three classical interpolation kernels via PSNR values.
Figure  \ref{Fig.sub.FirstN4D1 Better than cubic} and  Figure  \ref{Fig.sub.SecondtN4D1 Better than cubic} display the comparison cases  where the three  parameters of $S^4_{4/1}$ and $S^5_{4/1}$ are traversed, respectively.
The red points indicate cases where the PSNR values are higher than those of the three kinds of interpolation kernels.
Figure \ref{Fig.sub.FirstN4D1 Better than cubic}, \ref{Fig.sub.FirstN4D1 Better than cubic a03 thickness angle} and
Figure \ref{Fig.sub.SecondtN4D1 Better than cubic}, \ref{Fig.sub.SecondN4D1 Better than cubic a03 thickness angle}  show from two angles  that these points
corresponding  $S^4_{4/1}$ and $S^5_{4/1}$ fall roughly near a triangular region of a plane, respectively.
Fitting these red points in Figure \ref{Fig.sub.FirstN4D1 Better than cubic} results in a plane equation $a_{03}=-5.7392-2a_{01}-2.7906a_{02}$.
We selected two representative parameter sets in this plane: $(30,20,-121.5512)$ falls within the range of parameters we have explored, while $(80,100,-444.7992)$ lies outside.
Similarly, fitting these red points in Figure \ref{Fig.sub.SecondtN4D1 Better than cubic} leads to a plane equation $a_{03}=-5.7902-1.9574a_{01}-2.5645a_{02}$,
and we select two sets of parameters $(30,10,-90.1572)$ and $(50,10,-129.3052)$ in the plane.
The last four rows in Table
\ref{tab:best PSNR comparisions table For all kernels} show that
the PSNR values of these four sets of parameters are all greater than
those of the three classical interpolation kernels and Hu and Tan's interpolation kernel $S_{2/2}$ \cite{HuT06}.
The performance of $(80,100,-444.7992)$ may indicate that the optimal parameters of $S^4_{4/1}$ and $S^5_{4/1}$ are not limited to the range of parameters we have explored.


Next, we give   the experimental results for SSIM and FSIM.
Table \ref{tab:best SSIM comparisions table For all kernels}
and Table \ref{tab:best FSIM comparisions table For all kernels} demonstrate that optimal parameters $a_{02}$ of the kernel function $S_3(t)$ vary depending on the images used.
In contrast, our interpolation kernels using the four fixed sets of parameters remain competitive with the cubic interpolation kernel $S_3(t)$ in most cases.

The approximation order is a crucial criterion for evaluating interpolation kernels. However, real-world images are often complex and diverse, and the approximation order alone usually cannot directly measure the quality of image magnification.
As demonstrated in Table \ref{tab:best PSNR comparisions table For all kernels}, Table \ref{tab:best FSIM comparisions table For all kernels} and Table \ref{tab:best SSIM comparisions table For all kernels}, the parameter value of $a_{02}$  varies with different images when the cubic spline interpolation kernel $S_3$ achieves the highest values in terms of PSNR, SSIM, and FSIM values. Notably, it does not equal $-\frac{5}{2}$, failing to reach the corresponding highest approximation order.
In contrast, a rich family of functions may potentially provide more suitable interpolation kernels for complex real-world images. Although our interpolation kernel functions degenerate into polynomial forms when the approximation order exceeds $1$, they offer greater diversity and potential.
Our interpolation kernels $S^4_{4/1}$ and $S^5_{4/1}$, under the condition that the approximation order is 1 and with fixed parameters, outperform the best performances of the cubic interpolation kernel $S_3$  when the parameter $a_{02}$ varies within the interval $[-7, 1]$. In summary, the quartic/linear interpolation kernels $S^4_{4/1}$ and $S^5_{4/1}$ exhibit better performance than the cubic interpolation kernel $S_3(t)$.

\section{Conclusion}\label{sec:conc}
This paper introduces   parametric cubic/linear  $S_{3/1}(t)$ and  quartic/linear  $S_{4/1}(t)$  interpolation kernels supported on $[-2,2]$.
These interpolation kernels are all symmetric, $C^1$ continuous  and possess certain degrees of approximation order.
Numerical results demonstrate that one quartic/linear interpolation kernel can outperform the cubic interpolation kernel in terms of PSNR, SSIM and FSIM.
The simplicity and versatility of these kernels, along with their free parameters and favorable properties, make them a promising tool in approximation theories and CAGD.

Key direction of future work will be to construct interpolation kernel functions that not only possess a rational form but also achieve approximation orders exceeding the first order. This endeavor is expected to further enhance the flexibility and applicability of such kernels, paving the way for more robust and accurate applications in complex scenarios.

\section*{Data availability}
The Mathematica and MATLAB codes used for the theoretical derivations and numerical results in this paper are available  at  https://github.com/mathEuler/rationalPoly.
\section*{Acknowledgments}
This work was partly supported 
by the National Natural Science Foundation of China (grants 12261037) and 
the Scientific Research Fund of Liaoning Provincial Education Department (LQ2020020).

\section*{Declaration of generative AI and AI-assisted technologies in the writing process}
During the preparation of this work the authors used ChatGPT in order to improve language and readability. 
After using this tool, the authors reviewed and edited the content as needed and take full responsibility 
for the content of the publication.
\bibliographystyle{elsarticle-num}
\bibliography{rationalPolynomialImageInterpolation}

\end{document}